\documentclass[11pt]{amsart}

\usepackage[english]{babel}
\usepackage[T1]{fontenc}
\usepackage[utf8]{inputenc}

\usepackage{mathpazo}
\usepackage{graphicx}
\usepackage{microtype}
\usepackage{enumerate}
\usepackage{fontawesome}

\usepackage{csquotes}
\usepackage[dvipsnames]{xcolor}
\usepackage{url}

\usepackage{breakurl}
\usepackage[breaklinks,bookmarks, colorlinks, citecolor=MidnightBlue!75!black, urlcolor=MidnightBlue!75!black, linkcolor=MidnightBlue!75!black]{hyperref}

\calclayout

\AtBeginDocument{%
   \def\MR#1{}
}

\usepackage{scalerel}

\usepackage{amsmath, amssymb, amsthm, amscd}
\usepackage{commath, mathtools, mathrsfs}
\usepackage{thmtools}
\usepackage{tikz, tikz-cd}
\usetikzlibrary{matrix, arrows, decorations.pathmorphing}

\theoremstyle{plain}

\newtheorem{theorem}{Theorem}
\newtheorem{lemma}[theorem]{Lemma}
\newtheorem{proposition}[theorem]{Proposition}
\newtheorem{corollary}[theorem]{Corollary}

\newtheorem*{theorem*}{Theorem}
\newtheorem*{lemma*}{Lemma}
\newtheorem*{proposition*}{Proposition}
\newtheorem*{corollary*}{Corollary}

\theoremstyle{definition}

\newtheorem*{conjecture*}{Conjecture}
\newtheorem*{remark*}{Remark}
\newtheorem*{definition*}{Definition}
\newtheorem*{observation*}{Observation}

\newcommand{\on}{\operatorname}

\newcommand{\br}{\on{Br}}

\newcommand{\C}{\mathbb{C}}

\newcommand{\gal}{\on{Gal}}

\newcommand{\G}{\mathbb{G}}

\newcommand{\mc}{\mathcal}
\newcommand{\mf}{\mathfrak}

\newcommand{\pic}{\on{Pic}}

\newcommand{\Q}{\mathbb{Q}}

\newcommand{\R}{\mathbb{R}}

\newcommand{\spec}{\on{Spec}}

\newcommand{\upic}{\underline{\pic}}

\newcommand{\xar}{\xrightarrow}
\newcommand{\Z}{\mathbb{Z}}

\renewcommand{\epsilon}{\varepsilon}
\renewcommand{\H}{\on{H}}
\renewcommand{\hom}{\on{Hom}}

\renewcommand{\O}{\mc{O}}
\renewcommand{\P}{\mathbb{P}}
\renewcommand{\phi}{\varphi}

\renewcommand{\hat}{\widehat}

\usepackage[giveninits=true,style=alphabetic,backend=biber,maxnames=99,maxalphanames=5,giveninits=true]{biblatex}
\bibliography{qnf}

\usepackage{todonotes}

\title{On the section conjecture and Brauer-Severi varieties}
\date{}
\subjclass[2010]{}

\author{Giulio Bresciani}
\address{Freie Universit\"at Berlin, Arnimallee 3, 14195, Berlin, Germany}
\email{gbresciani@math.fu-berlin.de}
\thanks{The author is supported by the DFG Priority Program "Homotopy Theory and Algebraic Geometry" SPP 1786}

\begin{document}

\begin{abstract}
	J. Stix proved that a curve of positive genus over $\mathbb{Q}$ which maps to a non-trivial Brauer-Severi variety satisfies the section conjecture. We prove that, if $X$ is a curve of positive genus over a number field $k$ and the Weil restriction $R_{k/\mathbb{Q}}X$ admits a rational map to a non-trivial Brauer-Severi variety, then $X$ satisfies the section conjecture. As a consequence, if $X$ maps to a Brauer-Severi variety $P$ such that the corestriction $\operatorname{cor}_{k/\mathbb{Q}}([P])\in\operatorname{Br}(\mathbb{Q})$ is non-trivial, then $X$ satisfies the section conjecture.
\end{abstract}

\maketitle

Let $X$ be a geometrically connected variety over a field $k$ with separable closure $\bar{k}$, there is a short exact sequence of étale fundamental groups
\[0\to\pi_{1}(X_{\bar{k}})\to\pi_{1}(X)\to\gal(\bar{k}/k)\to 0.\]
Grothendieck's section conjecture predicts that, if $X$ is a smooth, projective curve of genus at least $2$ and $k$ is a number field, the set of rational points $X(k)$ is in natural bijection with sections of the sequence above modulo the action of $\pi_{1}(X_{\bar{k}})$ by conjugation.

Thanks to an idea of Tamagawa \cite{tam97} \cite[Corollary 102]{sti13} it is sufficient to prove the conjecture for curves with no rational points, and some results have been proved about such curves. The section conjecture holds for $X$ if

\begin{itemize}
	\item the number field $k$ has a real place $k\hookrightarrow \R$ such that $X_{\R}(\R)=\emptyset$, see \cite[Corollary 3.13]{moc03}, or
	\item the class of $\upic^{1}_{X}$ in $\H^{1}(k,\upic^{0}_{X})$ is not divisible, see \cite[Theorem 1.2]{hs09}, or
	\item $k=\Q$ and $X$ maps to a non-trivial Brauer-Severi variety, see \cite[Corollary 18]{sti10}.
\end{itemize}

The last condition, which is due to Stix, holds over any number field $k$, but with an additional hypothesis: for every prime number $p$, it is required that $X$ has bad reduction at most at one place $\mf{p}$ of $k$ over $p$, see \cite[Theorem 17]{sti10}. We provide a different generalization based on Weil's restriction of scalars.

Recall that, given a finite separable extension $k/h$ of fields and a quasi-projective variety $X$ over $k$, the Weil restriction $R_{k/h}X$ is a quasi-projective variety over $h$ characterized by a functorial bijection $\hom(S,R_{k/h}X)\simeq \hom(S_{k},X)$ for schemes $S$ over $h$. In particular, $k$-rational points of $X$ are in natural bijection with $h$-rational points of $R_{k/h}X$. 

If $X$ is a curve of genus $g$ over a number field $k$, then $(R_{k/\Q}X)_{\bar{\Q}}$ is a product of $[k:\Q]$ curves of genus $g$, so passing to the Weil restriction is basically a trade-off between the complexity of the base field and the complexity of the variety. We prove the following.

\begin{theorem}\label{main}
	Let $X$ be a smooth, projective, geometrically connected curve of positive genus over a number field $k$. Assume that $R_{k/\Q}X$ admits a rational map to a non-trivial Brauer-Severi variety. Then the section conjecture holds for $X$. Equivalently, if $\pi_{1}(X)\to\gal(\bar{k}/k)$ admits a section, then the map $\br(\Q)\to\br(R_{k/\Q}X)$ is injective.
\end{theorem}

As a consequence, we get the following corollaries.

\begin{corollary}\label{cor}
	Let $X$ be a smooth projective curve of positive genus over a number field $k$ and $P$ a Brauer-Severi variety such that the corestriction 
	$\on{cor}_{k/\Q}[P]\in\br(\Q)$ is non-trivial. If there exists a morphism $X\to P$, then the section conjecture holds for $X$.	Equivalently, if $\pi_{1}(X)\to\gal(\bar{k}/k)$ admits a section, the kernel of $\br(k)\to\br(X)$ is contained in the kernel of $\on{cor}_{k/\Q}:\br(k)\to\br(\Q)$. 
\end{corollary}

\begin{corollary}\label{res}
	Let $k$ be a number field and $P$ a Brauer-Severi variety over $\Q$ with $[k:\Q][P]\neq 0\in\br(\Q)$. If $X$ is a smooth projective curve over $\Q$ of positive genus with a morphism $X\to P$, then the section conjecture holds for the base change $X_{k}$.
\end{corollary}

Our argument is analogous to Stix's one and we rely heavily on his results. Our contribution consists essentially of two things: we realized that such a generalization was possible and we overcame the lack, for higher dimensional varieties, of a sufficiently strong analogue of Lichtenbaum's theorem about the period and index of a curve over a $p$-adic field, which is an essential ingredient of Stix's proof.

We mention that it is possible to prove \autoref{cor} (and thus \autoref{res}) analogously to Stix' theorem over $\Q$ \cite[Corollary 18]{sti10} without using Weil's restriction of scalars. The proof is basically the same plus the observation that, if $\alpha\in\br(k)$ is a Brauer class over the number field $k$, the Hasse invariant of $\on{cor}_{k/\Q}(\alpha)$ at $p$ is the sum of the Hasse invariants of $\alpha$ at places over $p$.

\section{Weil restriction and the section conjecture}

The behaviour of the étale fundamental group and the section conjecture with respect to the Weil restriction of scalars has been studied by J. Stix in \cite{sti10b}. Let $k/h$ be a finite, separable extension of fields and $X$ is a geometrically connected variety over $k$. Assume either that $X$ is proper or that $k$ has characteristic $0$. Stix describes explicitly the étale fundamental group of $R_{k/h}X$ in terms of the one of $X$, and uses this description to show that the section conjecture holds for $X$ if and only if it holds for $R_{k/h}X$. We give here an alternative treatment based on \emph{étale fundamental gerbes}.

Recall that A. Vistoli and N. Borne have introduced the étale fundamental gerbe $X\to\Pi_{X/k}$ of a geometrically connected scheme, see \cite[Section 9]{bv15} and \cite[Appendix]{bre21}. The set of Galois sections of the étale fundamental group is in natural bijection with the isomorphism classes of $\Pi_{X/k}(k)$. We show that the étale fundamental gerbe and the Weil restriction commute.

\begin{proposition}\label{weilfun}
	Let $k/h$ be a finite separable extension of fields, and $X$ a geometrically connected quasi-projective variety over $k$. Assume either that $X$ is proper or that $\on{char}k=0$.
	
	Then $R_{k/h}X$ is geometrically connected and the natural morphism $R_{k/h}X\to R_{k/h}\Pi_{X/k}$ induces a natural isomorphism
	\[\Pi_{R_{k/h}X/h}\xar{\sim}R_{k/h}\Pi_{X/k}.\]
	\begin{proof}
		Let $\bar{h}/h$ be a separable closure. We have natural isomorphisms $(R_{k/h}X)_{\bar{h}}=\prod\sigma^{*}X$ and $(R_{k/h}\Pi_{X/k})_{\bar{h}}\simeq\prod\sigma^{*}\Pi_{X/k}=\prod\Pi_{\sigma^{*}X/k}$ where the product runs over the $h$-linear embeddings $\sigma:k\subset \bar{h}$, see \cite[Theorem 1.3.2]{wei82} (Weil's original work deals only with varieties, but his proof easily generalizes to any fibered category). In particular $R_{k/h}X$ is geometrically connected and $R_{k/h}\Pi_{X/k}$ is a pro-finite étale gerbe, thus the morphism $R_{k/h}X\to R_{k/h}\Pi_{X/k}$ induces a natural morphism
		\[\phi:\Pi_{R_{k/h}X/h}\to R_{k/h}\Pi_{X/k}.\]
		The base change of $\phi$ to $\bar{h}$ is an isomorphism since both terms are naturally isomorphic to $\prod_{\sigma}\Pi_{\sigma^{*}X/k}$, thus $\phi$ is an isomorphism too. 
	\end{proof}
\end{proposition}

\begin{corollary}{\cite[Theorem 2, Theorem 3]{sti10b}}
	The set of isomorphism classes of Galois sections of $\pi_{1}(X)$ over $k$ is in natural bijection with the one of $\pi_{1}(R_{k/h}X)$ over $h$.\qed
\end{corollary}

\begin{corollary}{\cite[Theorem 4]{sti10b}}
	Let $X$ be a smooth, projective curve of genus $g\ge 2$ over a number field $k$. The section conjecture holds for $X$ if and only if it holds for $R_{k/\Q}X$.\qed
\end{corollary}

If $s\in\Pi_{X/k}(k)$ is a Galois section, denote by $s^{\Q}\in \Pi_{R_{k/\Q}X/\Q}(\Q)=R_{k/\Q}\Pi_{X/\Q}(\Q)$ the induced section. Recall that an étale neighbourhood of $s$ is a finite étale cover $Y\to X$ such that $s$ lifts to $\Pi_{Y/k}(k)$. 

\begin{corollary}\label{weilnb}
	Let $X$ be a smooth, projective curve of positive genus over a number field $k$ and $s\in\Pi_{X/k}(k)$ a Galois section. If $Y\to X$ is an étale neighbourhood of $s$, then $R_{k/\Q}Y\to R_{k/\Q}X$ is an étale neighbourhood of $s^{\Q}$. The étale neighbourhoods of this form are cofinal in the system of all étale neighbourhoods of $s^{\Q}$.\qed
\end{corollary}

\section{Morphisms to Brauer-Severi varieties}

If $X$ is a scheme over $k$, denote by $\br(X/k)$ the kernel of $\br(k)\to\br(X)$. If $X$ is a regular variety, the restriction map $\br(X)\to\br(k(X))$ is injective \cite[Corollary IV.2.6]{mil80} and thus $\br(X/k)=\br(k(X)/k)$. In particular, a Brauer class $[P]\in\br(k)$ of a Brauer-Severi variety $P$ is in $\br(X/k)$ if and only if there exists a rational map $X\dashrightarrow P$.

If $X$ is a smooth, projective variety, the Leray spectral sequence in étale cohomology for the map $X\to\spec k$ gives a short exact sequence
\[0\to\pic(X)\to\upic_{X}(k)\to\br(X/k)\to 0\]
where $\upic_{X}$ is the Picard scheme of $X$ and $\pic(X)=\H^{1}(X,\G_{m})$ is the Picard group. Let us call $\beta$ the homomorphism $\upic_{X}(k)\to\br(X/k)$.

\begin{lemma}\label{brdiv}
	Let $X$ be a smooth, projective variety over a field $k$ of characteristic $0$, $s\in\Pi_{X/k}(k)$ a section, $b\in\br(X/k)$ a Brauer class split by $X$. Assume that the second étale homotopy group of $X_{\bar{k}}$ is trivial. For every positive integer $n$, there exists an étale neighbourhood $f:Y\to X$ of $s$ such that $f^{*}b\in n\br(Y/k)$.
	\begin{proof}
		Let $L\in\upic_{X}(k)$ be such that $\beta(L)=b\in\br(X/k)$, and $L_{\bar{k}}\in\pic(X_{\bar{k}})$ the associated line bundle over $X_{\bar{k}}$. We have an exact sequence
		\[\H^{1}(X_{\bar{k}},\mu_{n})\to\pic(X_{\bar{k}})\xar{\cdot n}\pic(X_{\bar{k}})\xar{\delta}\H^{2}(X_{\bar{k}},\mu_{n}).\]
		Let $\on{\acute{E}t}(X_{\bar{k}})$ be the étale homotopy type of $X_{\bar{k}}$ and $\on{cosk}_{3}(\on{\acute{E}t}(X_{\bar{k}}))$ its third coskeleton, since $\pi_{2}^{\rm\acute{e}t}(X_{\bar{k}})$ is trivial we have $\on{cosk}_{3}(\on{\acute{E}t}(X_{\bar{k}}))=K(\pi_{1}^{\rm\acute{e}t}(X_{\bar{k}}),1)$. Therefore, we have
		\[\H^{2}(X_{\bar{k}},\mu_{n})=\H^{2}(\on{\acute{E}t}(X_{\bar{k}}),\mu_{n})=\H^{2}(\on{cosk}_{3}(\on{\acute{E}t}(X_{\bar{k}})),\mu_{n})=\H^{2}(\pi_{1}^{\rm\acute{e}t}(X_{\bar{k}}),\mu_{n}),\]
		see \cite[Corollary 9.3]{am69} for the first equality. Since the base change to $\bar{k}$ of the étale neighbourhoods of $s$ are cofinal in all finite étale covers of $X_{\bar{k}}$, there exists an étale neighbourhood $g:X'\to X$ of $s$ such that $g_{\bar{k}}^{*}\delta(L_{\bar{k}})=0$ and thus $g_{\bar{k}}^{*}L_{\bar{k}}\in\pic(X'_{\bar{k}})$ is divisible by $n$.
		
		Choose $M\in\pic(X'_{\bar{k}})$ such that $M^{n}=g_{\bar{k}}^{*}L_{\bar{k}}$. Since the Picard scheme of $X'$ is locally of finite type, the residue field of $\spec\bar{k}\xar{M}\upic_{X'}$ is finite over $k$ and thus the Galois orbit of $M$ is finite. If $\sigma\in\gal(\bar{k}/k)$ is an element, since $g_{\bar{k}}^{*}L_{\bar{k}}$ is Galois invariant then $M\otimes \sigma M^{-1}$ is $n$-torsion, and thus it comes from $\H^{1}(X'_{\bar{k}},\mu_{n})=\hom(\pi_{1}^{\rm\acute{e}t}(X'_{\bar{k}}),\mu_{n})$. It follows that there exists an étale neighbourhood $h:Y\to X'$ of $s$ such that $h_{\bar{k}}^{*}M\in\pic(Y_{\bar{k}})$ is Galois-invariant and thus descends to an element $N\in\upic_{Y}(k)$.
		
		Let $f:Y\to X$ be the composition, we have $N^{n}=f^{*}L$ and hence $n\beta(N)=\beta(f^{*}L)=f^{*}b$.
	\end{proof}
\end{lemma}

Recall that a group $G$ is \emph{good in the sense of Serre} if $\H^{i}(\hat{G},M)\to\H^{i}(G,M)$ is an isomorphism for every finite $G$-module $M$, see \cite[I.2.6]{ser94}. Fundamental groups of complex curves are good \cite[Proposition 3.6]{gjz08}.

\begin{lemma}\label{topet}
	Let $X$ be a variety over $\C$. Assume that $\pi_{2}^{\rm top}(X^{\rm an})$ is trivial and that $\pi_{1}^{\rm top}(X^{\rm an})$ is good in the sense of Serre. Then $\pi_{2}^{\rm \acute{e}t}(X)$ is trivial.
	\begin{proof}
		Write $\pi_{i}=\pi_{i}^{\rm\acute{e}t}(X)$. The hypothesis implies that the natural homomorphism $\H^{2}(\pi_{1},M)\to\H^{2}(X,M)$ is bijective for every finite $\pi_{1}$-module $M$.
		
		Assume by contradiction that $\pi_{2}$ is not trivial, then there exists a finite homotopy type $F$ with $\pi_{n}(F)=0$ for $n\ge 3$ and a map $\on{\acute{E}t}(X)\to F$ such that $\pi_{2}\to\pi_{2}(F)$ is non-trivial. Up to passing to finite étale coverings of $X$ and $F$, we may assume that $\pi_{1}(F)$ is trivial and hence $F=K(M,2)$ for some finite abelian group $M$ (the fundamental group of the covering of $X$ is still good thanks to \cite[Lemma 3.2]{gjz08}).
		 
		We thus have a map $X\to K(M,2)$ inducing a non-trivial homomorphism $\pi_{2}\to M$. This defines a cohomology class $\alpha\in\H^{2}(X,M)$ not in the image of $\H^{2}(\pi_{1},M)\to\H^{2}(X,M)$, and this is absurd.
	\end{proof}
\end{lemma}

\begin{lemma}\label{weilcor}
	Let $k/h$ be a finite separable extension and $P/k$ a Brauer-Severi variety. There exists a Brauer-Severi variety $Q/h$ with $[Q]=\on{cor}_{k/h}([P])\in\br(h)$ and a closed embedding $R_{k/h}P\hookrightarrow Q$.
	\begin{proof}
		Let $\bar{h}/h$ be a separable closure, then
		\[R_{k/h}P=\left(\prod_{\sigma}\sigma^{*}P\right)/\gal(\bar{h}/h)\]
		where the product runs over $h$-linear embeddings $\sigma:k\to \bar{h}$. The Galois group $\gal(\bar{h}/h)$ permutes the factors and the stabilizer $\gal(\bar{h}/\sigma k)$ of $\sigma^{*}P$ acts on it. Note that, even though $\sigma^{*}P$ is a projective space over $\bar{h}$, the action of $\gal(\bar{h}/\sigma k)$ is non-standard.
		
		The external tensor product $\boxtimes_{\sigma}\O(1)\in\pic(\prod\sigma^{*}P)$ is naturally $\gal(\bar{h}/h)$-equivariant and thus the Segre embedding
		\[S:\prod_{\sigma}\sigma^{*}P\hookrightarrow \P(\H^{0}(\boxtimes_{\sigma}\O(1)))\]
		is naturally $\gal(\bar{h}/h)$-equivariant.	The quotient $Q=\P(\H^{0}(\boxtimes_{\sigma}\O(1)))/\gal(\bar{h}/h)$ is a Brauer-Severi variety over $h$.  
		
		Using the fact that summation in the Brauer group can be computed using the Segre embedding of Brauer-Severi varieties \cite[\S 4]{art82} and the fact that the corestriction homomorphism is the derived augmentation homomorphism, it is easy to show that that $[Q]=\on{cor}_{k/h}([P])$. Moreover, the Segre embedding $S$ descends to a closed embedding $R_{k/h}P\hookrightarrow Q$ since it is $\gal(\bar{h}/h)$-equivariant.
	\end{proof}
\end{lemma}

\section{Proof of the main theorem}

Let us now prove \autoref{main}. Let $X$ be a smooth projective curve over a number field $k$ such that $R_{k/\Q}X$ admits a rational map to a non-trivial Brauer-Severi variety, we want to show that $\Pi_{X/k}(k)$ is empty. Assume by contradiction that there exists a section $s\in\Pi_{X/k}(k)$ and let $b\in\br(R_{k/\Q}X/\Q)$ be a non-trivial Brauer class. Let $s^{\Q}\in\Pi_{R_{k/\Q}X/\Q}(\Q)$ be the associated section.

Since $(R_{k/\Q}X)_{\bar{k}}$ is a product of curves and the fundamental group of a curve over $\bar{k}$ is good in the sense of Serre \cite[Proposition 3.6]{gjz08}, then \autoref{topet} implies that $\pi_{2}^{\rm \acute{e}t}((R_{k/\Q}X)_{\bar{k}})$ is trivial and we may thus apply \autoref{brdiv} to $R_{k/\Q}X$ and $s^{\Q}$. If we apply \autoref{brdiv} together with \autoref{weilnb}, for every $N>0$ we may find an étale neighbourhood $X_{N}\to X$ of $s$ and a Brauer class $b_{N}\in\br(R_{k/\Q}X_{N}/\Q)$ such that $Nb_{N}=b\in\br(\Q)$ is non-trivial. 

Let $l/k/\Q$ be a Galois closure. Up to replacing $X$ with $X_{2[l:\Q]}$ and $b$ with $b_{2[l:\Q]}$, we may assume that $2[l:\Q]b\in\br(\Q)$ is non-trivial.

Fix $p$ a prime number, let us show that the order of the Brauer class $2[l:\Q]b_{\Q_{p}}$ is a power of $p$. Let $L$ be the completion of $l$ at some place over $p$, we have that $L/\Q_{p}$ is a Galois extension such that $[L:\Q_{p}]$ divides $[l:\Q]$, it is enough to show that the order of $[L:\Q_{p}]b_{\Q_{p}}$ is a power of $p$. Let $\Sigma$ be the set of embeddings $k\to L$, we have
\[(R_{k/\Q}X)_{L}=\prod_{\sigma\in\Sigma}\sigma^{*}X.\]
The section $s\in\Pi_{X/k}(k)$ induces a section $\sigma^{*}s\in\Pi_{\sigma^{*}X/L}(L)$ for every embedding $\sigma:k\to L$. By \cite[Theorem 15]{sti13}, this implies that the index of $\sigma^{*}X$ is a power of $p$ for every $\sigma$. Let $D_{\sigma}\in Z_{0}(X_{\sigma})$ be a $0$-cycle whose degree is a power of $p$, then $\otimes_{\sigma} D_{\sigma}$ is a $0$-cycle on $\prod_{\sigma}\sigma^{*}X$ whose degree is a power of $p$. It follows that the index of $\prod_{\sigma}\sigma^{*}X$ is a power of $p$, too.

Since $\Q(R_{k/\Q}X)$ splits $b$, there exists a Brauer-Severi variety $P$ with $[P]=b$ and a smooth projective variety $Y/\Q$ birational to $R_{k/\Q}X$ with a morphism $Y\to P$. Since the index is a birational invariant, the index of $Y_{L}$ is a power of $p$, it follows that the index of $b_{L}=[P_{L}]$ is a power of $p$. This implies that the order (i.e. the period) of $b_{L}\in\br(L)$ is a power of $p$, and finally the same holds for $\on{cor}_{L/\Q_{p}}b_{L}=[L:\Q_{p}]b_{\Q_{p}}$.

We thus have that the order of $2[l:\Q]b_{\Q_{p}}$ is $p$-primary for every $p$, and clearly $2[l:\Q]b_{\R}=0\in\br(\R)=\Z/2\Z$. 

The rest of the argument is analogous to Stix's one. Let $\alpha_{p}\in\Q/\Z$ be the Hasse invariant of $2[l:\Q]b_{\Q_{p}}$, by the Brauer-Hasse-Noether theorem we have $\sum_{p}\alpha_{p}=0\in\Q/\Z$. Since $\alpha_{p}$ is $p$-primary for every $p$, it follows that $\alpha_{p}=0$ for every $p$ and thus $2[l:\Q]b\in\br(\Q)$ is trivial, which is absurd. This concludes the proof of \autoref{main}.

\subsection*{Corollaries} \autoref{main} implies \autoref{cor} using \autoref{weilcor}. Since the composition $\on{cor}_{k/\Q}\circ\on{res}_{k/\Q}:\br(\Q)\to\br(\Q)$ is multiplication by $[k:\Q]$, \autoref{cor} implies \autoref{res}.

\subsection*{Acknowledgements}

\autoref{brdiv} was found during joint work with A. Vistoli. I would like to thank an anonymous referee for many useful remarks.

\printbibliography
	
\end{document}